\documentclass[a4paper, 12pt]{article}
\usepackage{tgheros}
\usepackage{amsthm}
\newtheorem{thm}{Theorem}
\newtheorem{cor}[thm]{Corollary}

\newtheorem{lemma}[thm]{Lemma}
\newtheorem{defn}{Definition}

\usepackage{amssymb,amsmath}
\usepackage[all]{xy}

\newtheorem*{theorem*}{Theorem}

\DeclareMathOperator{\F}{\mathbb{F}}

\begin{document}
\baselineskip=16.3pt
\parskip=14pt

\begin{center}
\section*{Linearization of polynomials in prime characteristic, with applications to the Golay code and Steiner system}

{\large 
Rod Gow and Gary McGuire
 \\ { \ }\\
School of Mathematics and Statistics\\
University College Dublin\\
Ireland}
\end{center}

 \subsection*{Abstract}
 
Let $F$ be any field containing the finite field of order $q$. A $q$-polynomial $L$ over $F$ is an element of the polynomial ring $F[x]$ with the property that all powers of $x$ that appear in $L$ with nonzero coefficient have exponent a power of $q$. It is well known that given any ordinary polynomial $f$ in $F[x]$, there exists a $q$-polynomial  that is divisible by $f$. We study the smallest degree of such a $q$-polynomial. This is equivalent to studying the $\F_q$-span of the roots of $f$ in a splitting field.  We relate this quantity to the representation theory of the Galois group of $f$.  As an application we give a simultaneous construction of the binary Golay code of length 24,  and the Steiner system on 24 points.
    
  
  
  \vskip1in
  
  {\bf Keywords} Galois group, linearized polynomial

  \newpage
  
\section{Introduction}
  
\noindent Let $p$ be a prime number, and let $q$ be a positive power of $p$.
Let $\F_q$ denote the finite field with $q$ elements.
Let $F$ be a field of characteristic $p$,
and assume that $F$ contains $\F_q$.
A $q$-linearized polynomial over $F$  is a polynomial of the form 
\begin{equation}\label{lin1}
L=a_nx^{q^n}+a_{n-1}x^{q^{n-1}}+\cdots+a_2x^{q^2}+a_1x^q+a_0 x \in F[x].
\end{equation}
If $a_n\not=0$ we say that $n$ is the $q$-degree of $f$.
We will usually say that $L$ is $q$-polynomial instead of a $q$-linearized polynomial.

The set of roots of $L$ forms an $\F_q$-vector space,
which is contained in a splitting field of $L$. 
We make this statement more precise in the following lemma.

\begin{lemma} \label{vector_space}
Let $F$ be a field of prime characteristic $p$ that contains $\F_q$. Let $L$ be a $q$-polynomial of $q$-degree $n$ in $F[x]$, with
\[
L=a_nx^{q^n}+a_{n-1}x^{q^{n-1}}+\cdots+a_2x^{q^2}+a_1x^q+a_0 x.
\]
Let $E$ be a splitting field for $L$ over $F$ and let $V$ be the set of roots of $L$ in $E$. Let $G$ be the Galois group of $E$
over $F$. Suppose that $a_0\neq 0$. Then $V$ is an $\F_q$-vector space of dimension $n$ and $G$ is naturally 
isomorphic to a subgroup of $GL(n,q)$.
\end{lemma}

Let $f\in F[x]$ be an arbitrary non-constant polynomial, and suppose that $f$ divides a
$q$-polynomial $L$ (it is well known that this always happens, see Lemma \ref{mind}).
Abhyankar and Yie \cite{AY} 
say that $f\in F[x]$ linearizes at $n$ if there exists a 
$q$-linearized polynomial
$L$ of $q$-degree $n$ such that $f$ divides $L$.
If $f$ linearizes at $n$ where $n$ is significantly smaller than $\deg(f)$, then useful information 
can be obtained about the Galois group of $f$, as it may then be considered as a subgroup
of $GL(n,q)$.
They refer to this as Serre's linearization trick. We discuss this idea in Section \ref{SLT}.

A link to representation theory also follows from Lemma \ref{vector_space}. 
Given $f$ in $F[x]$ with Galois group $G$, 
if there exists a 
$q$-linearized polynomial
$L$ of $q$-degree $n$ such that $f$ divides $L$,
then $G$ has a representation over $\F_q$ of dimension $n$.
We will explain in Section \ref{perm_module_section} how interesting representations for $G$ arise
from polynomials that linearize at $n$ where $n$ is significantly smaller than $\deg(f)$.
We explain this starting with the permutation module for $G$.
Section \ref{possible_min} considers which representations of $G$ over $\F_q$ 
can appear in this way.

In Section \ref{sect_sm} 
we briefly consider the minimal linearization for polynomials whose
Galois group is the symmetric group. This is the generic situation.

Section \ref{golay} contains our  construction of both the Golay code
and the Steiner system, on 24 points.
This code and design is the classic example of an error-correcting 
code and a combinatorial design being linked.
The Steiner system on 24 points was discovered first by Carmichael \cite{C1} in 1931,
and the Golay code was first discovered by Golay in 1949.
Carmichael shows that the Steiner system arises from the action of 
$M_{24}$ on 24 points.
That the Golay code, Mathieu group and Steiner system (on 23 points)
are related was first written about by Paige \cite{P}, and further results were given
  in the 1960s by Assmus and Mattson \cite{AM}.
There are now many different ways  to construct the Golay code, and the Steiner system.
In Section \ref{golay} we provide one more! We use Galois theory over $\F_2(t)$ to
construct both objects together, and the arguments are elementary and natural in the context
of the splitting field of a polynomial.  Our proof is self-contained and does not assume knowledge of $M_{24}$.

In Section \ref{m24} we relate our work to known facts about the Mathieu group $M_{24}$
and its representations over $\F_2$. This is not self-contained, and provides context.
We prove that there exists a polynomial with minimal linearization of 2-degree 23.

Section 8 briefly outlines a similar approach for the affine general linear group and the 
Hamming code of length 8, and Section 9 discusses the projective linear group $PSL(2,p)$.
These last two sections are brief, but they indicate that the approach does work for
other groups.

\section{Background and Serre's Linearization Trick}\label{SLT}

\noindent We present a simple result which we will need later, it is a converse to Lemma \ref{vector_space}.
The proof is essentially the same as Theorem 3.52  in \cite{LN},
which goes back to Dickson.

\begin{lemma}\label{lin}
Let $V$ be a finite dimensional vector space over $\F_q$, which is contained in a field extension 
$E$ of $\F_q$.
Then the polynomial $\prod_{v\in V} (x-v)$ is a $q$-linearized polynomial. 
\end{lemma}

\begin{proof}
Let $\alpha_1, \ldots, \alpha_n \in E$ be a basis for $V$.
Consider the polynomial in $E[x]$
\[
D(x):=\det
 \left[ \begin{array}{ccccc}
   \alpha_1&\alpha_1^q&\alpha_1^{q^2}&\cdots&\alpha_1^{q^{n}}\\
 \alpha_2&\alpha_2^q&\alpha_2^{q^2}&\cdots&\alpha_2^{q^{n}}\\
 \vdots&\vdots&\vdots&\ddots&\vdots\\
  \alpha_n&\alpha_n^q&\alpha_n^{q^2}&\cdots&\alpha_n^{q^{n}}\\
  x&x^q&x^{q^2}&\cdots&x^{q^{n}}\\
 \end{array} \right]
\]
which  will clearly be a $q$-polynomial in $x$, of $q$-degree at most $n$.
We claim that the roots of $D(x)$ are precisely the elements of $V$,
 from which it follows that $D(x)$ is a scalar multiple of  $\prod_{v\in V} (x-v)$.
 
 The claim follows by observing that each $\alpha_i$ is a root of $D(x)$,
 and since $D(x)$ is a $q$-polynomial, all $\F_q$-linear combinations of
 the $\alpha_i$  are roots of $D(x)$ by Lemma  \ref{vector_space}.
 Since $V$ has $q^n$ elements, and $D(x)$ has degree at most $q^n$,
 it follows that $D(x)$ has degree exactly $q^n$ and the proof is complete.
\end{proof}

Given any polynomial $f\in F[x]$, it is well known that we can find a $q$-polynomial
in $F[x]$ that is divisible by $f$.
This result goes back to Ore. Proofs can be found in Goss \cite{Goss} and in
Lidl and Niederreiter  \cite{LN} (page 120).
Here is the precise statement.

\begin{lemma}\label{mind}
Let $F$ be a field of characteristic $p$ that contains $\F_q$.
For any  polynomial $f\in F[x]$ of degree $m$,
there exists  $L\in F[x]$  which is divisible by $f$ and has  the form
\[
L(x)=\sum_{i=0}^d b_i x^{q^i}
\]
where $d\le m$. 
\end{lemma}

If the polynomial $f$ in Lemma \ref{mind} has repeated roots, then the same must be true of the $q$-polynomial
$L$, and this implies that the coefficient of $x$ in $L$ is 0. Conversely, as shall show shortly,
if $f$ has no repeated roots, the same is true of $L$ and hence the coefficient of $x$ in $L$
is nonzero. It is usual to restrict to the case that $f$ has distinct roots in this type
of work.

\begin{defn}
Let $f$ be a polynomial in $F[x]$ with no repeated roots. We let $m(f)$ denote the smallest
$q$-degree of a nonzero monic $q$-polynomial in $F[x]$ that is divisible by $f$.
\end{defn}

Lemma \ref{mind} implies that $m(f)\le \deg(f)$.

\begin{lemma}\label{mind2}
Let $F$ be a field of characteristic $p$ that contains $\F_q$.
Let $f$ be a polynomial in $F[x]$ with no repeated roots. Then
\begin{enumerate}
\item
The integer $m(f)$  is 
equal to the dimension of the $\F_q$-span of the roots of $f$ in a splitting field. 
\item If
$L$ is a monic $q$-polynomial of $q$-degree $m(f)$,
and if $L$ is divisible by $f$, then $L$ divides any other
$q$-polynomial in $F[x]$ that is  divisible by $f$. 
 \item The coefficient of $x$ in $L$ is nonzero.
 \end{enumerate}
\end{lemma} 

\begin{proof}
Let $E$ be a splitting field for $f$ over $F$.  The set of all 
$\F_q$-linear combinations of the roots of $f$  is
an $\F_q$-vector space $V$, say,  contained in $E$. 
Let $L$ be the monic polynomial $\prod_{v\in V} (x-v)$.

Then $L$ is a $q$-polynomial  by Lemma \ref{lin}, and $f$ 
divides $L$ because its roots are a subset of the roots of $L$ (here we apply the hypothesis
that $f$ has no repeated roots).
The $q$-degree of $L$ is  $\dim V$. 

Let $L'$ be any other nonzero $q$-polynomial in $F[x]$ divisible by $f$ and let $V'$ be its space
of roots in a splitting field over $F$. Since $f$ divides $L'$, the roots of $f$ are contained in $V'$. Thus since $V$ is spanned by the roots of $f$, $V$ is an $\F_q$-subspace of $V'$. 
We set $M=\prod_{v\in V'} (x-v)$. By previous arguments, $M$ is a $q$-polynomial
over $F$ divisible by $L$. Furthermore, $M$ divides $L'$, as its roots
are a subset of those of $L'$ (the roots of $L'$ may be repeated, but the different roots of $L'$ are those of $M$).

Given this, it is clear that $L$ has degree $m(f)$ and it divides any other
$q$-polynomial in $F[x]$ divisible by $f$. The construction of $L$ implies that
$m(f)$ equals the dimension of the $\F_q$-span of the roots of $f$.

$L$ must have distinct roots, because the coefficient of $x$
is the product of all the nonzero elements of $V$, and is therefore nonzero.
\end{proof}

\begin{defn}
Let $f$ be a polynomial in $F[x]$ with no repeated roots. We call the unique monic
$q$-polynomial of $q$-degree $m(f)$ that is divisible by $f$ 
the minimal $q$-polynomial of $f$ over $F$.
\end{defn}

The minimal $q$-polynomial $L$, say, of $f$ may have several irreducible factors
but cannot be irreducible, because $x$ is a factor. It can happen
that $L(x)/x$ is irreducible and that $L(x)/x$ is the original $f$ that we chose.

Here are some trivial bounds on the degree of the minimal $q$-polynomial, and yet
they are tight.

 \begin{lemma}\label{bounds}
 Let $f$ be a polynomial in $F[x]$ with $\deg f=m$. Suppose that $x$ does not divide $f$ 
 and $f$ has no repeated roots.
 Let $m(f)$ be the $q$-degree of its minimal $q$-polynomial $L$ over $F$. Then we have
 $\log_q (m+1) \le m(f) \le m$.
 Furthermore, these bounds are best possible.
 \end{lemma}
 
 \begin{proof}
 That $m(f)\le m$ follows from Lemma \ref{mind}. 
 We will prove in Theorem  \ref{the_value_of_m(f)} that equality can hold
 (and often does)  in the inequality $m(f)\le m$.

 For the other inequality, the proof of Lemma \ref{mind2} shows that
 $L$ has exactly $q^{m(f)}-1$ different nonzero roots. Since the roots of $f$ are a subset
 of those $L$, and we are assuming that they are different and nonzero, $m\le q^{m(f)}-1$. 
 This gives the required inequality.

 Let $M$ be a $q$-polynomial, and let
  $f(x)=M(x)/x$. 
 Let $\deg f=m$. 
 Then the minimal $q$-polynomial of $f$ is $M$,
 and this shows that equality can hold in the inequality $m\le q^{m(f)}-1$.
 \end{proof}

 We conclude this section with an observation about the Galois group of
 the irreducible factors of the minimal $q$-polynomial of an irreducible $f$.
 
\begin{lemma}\label{other_gg}
Let $f$ be an irreducible polynomial in $F[x]$ of degree
$m\ge 3$ with Galois group $G$.
Let $L$ be the minimal $q$-polynomial of $f$ over $F$.
Let $h\not= f$ be an irreducible factor of $L$ over $F$.
Then  the Galois group of $h$  is a homomorphic image of $G$.
\end{lemma}

\begin{proof}
The roots of the irreducible factors of $L$ are, by Lemma \ref{mind2},
$\F_q$-linear combinations of the roots of $f$.
The splitting field of an irreducible  factor $h$ is therefore a subfield of the splitting field of $f$.
By the Galois correspondence, $Gal(h)$ is a quotient of $Gal(f)$.
\end{proof}

Sometimes  the irreducible factors $h\not= f$ of $L$ have the same splitting field 
as $f$ and the same Galois group as $f$.
One instance of this occurs when the Galois group of $f$ is a simple group.
However, the general problem of characterizing homomorphic images 
 is a difficult one, which indicates that completely characterizing
the $m(f)$ values and minimal $q$-polynomials is a difficult problem.

\begin{table}
\begin{center}
\begin{tabular}{|c|c|c|c|c|c|c|}
\hline
\bf{Polynomial} &  $q$ & $m(f)$ & \bf{Galois Group} &  \bf{Minimal Linearization} \\
\hline
$x^{2}-t$	&  odd &$1$& $C_2$   & $x^q-t^{(q-1)/2}x$  \\
\hline
$x^4+(t+1)x^2+1$	&  3  &2 &   $C_2\times C_2$  &$x^{9}+(t^3+t-1)x^3+(t^2-t)x$  \\
\hline
$x^q+x+t$	&  any  &2 &   order $2q$   &$x^{q^2}+(1-t^{q-1})x^q-t^{q-1}x$  \\
\hline
$x^q+tx+1$	&  any  &2 &  order $q(q-1)$  &$x^{q^2}+(t^q-1)x^q-tx$  \\
\hline
$(x^q-x)^{q-1}+t$	&  any  &2 &  $AGL(1,q)$  &$x^{q^2}+(t-1)x^q-tx$  \\
\hline
$x^{q+1}+x+t$	&  any &$q$& $PGL(2,q)$   & too large \\
\hline
$x^{q+1}+tx+1$	&  any &$q$& $PSL(2,q)$   &too large  \\
\hline
$x^{23}+x^3+t$	&  $2$ & 11 & $M_{23}$   & 
$x^{2048} + t^{64}x^{256} + x^{128} + t^{32}x^{32} + t^{80}x^8 +$\\
\hline
&&&&$ t^{88}x^4 + t^{82}x^2 +t^{89}x$\\
\hline
$f_1 $	&  $2$ & 11 & $M_{24}$   & 
$x^{2048} + t^{64}x^{512} + t^8x^{16} + t^{16}x^8 + x  $\\
\hline
$x^{24}+x+t$	&  $2$ & 12 & $M_{24}$   & 
$x^{4096} + (t^{24} + t)x^{2048} + t^{128}x^{1024} +  $\\
\hline
&&&&$ (t^{88} +t^{65})x^{512} + t^{16}x^{32} +t^9x^{16} + $\\
\hline
&&&&$(t^{40} + t^{17})x^8 + x^2 + (t^{24} + t)x$\\
\hline
$f_2 $	&  $3$ & 5 & $M_{11}$   & 
$  x^{243} -t^{18}x^{27} + t^3x^9 -t^9x^3 -x$\\
\hline
$x^{3}+x+t$	&  7 &$2$& $S_3$   &
$x^{49} +  (-t^{14} -t^{10} + 3t^8 -t^6 + 1)x^7  +$ \\
\hline
&&&&$(t^{12} + 2t^{10} -t^8 -t^6)x$\\
\hline
\end{tabular}
\caption{\label{tab:experiments1} Polynomials with base field $\F_q (t)$, some calculated using magma. Here $f_1=x^{276} + t^8x^{84} + x^{69} + x^{46} + 1$ and
$f_2=x^{110} + t^3 x^{74} + 2x^{66} + x^{44} + 2t^6 x^{38} + 2t^3 x^{30} + 2x^{22} + t^3 x^8 + 2t^9 x^2 + 2$
}
\end{center}
\end{table}

Some polynomials $f\in \F_q(t)[x]$ have minimal $q$-polynomial of a nice form,
and we list some examples in Table \ref{tab:experiments1}.
However, for many $f$ the minimal $q$-polynomial does not seem to exhibit any patterns.
The minimal 2-polynomial for $x^{23}+x+t $ over $\F_2(t)$  is enormous, whereas the 
minimal 2-polynomial for $x^{24}+x+t$ has very few terms and is in  Table \ref{tab:experiments1}.
Even an innocuous looking polynomial like $x^3+x+t$ does not have a predictable minimal $q$-polynomial;
we give an example when $q=7$ in Table \ref{tab:experiments1}.
Lemma \ref{anyLpoly}  shows in many case that a minimal $q$-polynomial over $\F_q(t)$ can  
be arbitrarily complicated.

 \section{Permutation Modules}\label{perm_module_section}

\noindent Let $f\in F[x]$ have degree $m$ and no repeated roots.
 Let $E$ be a splitting field for $f$ over $F$ and let $G$ be the Galois group of $f$ over $F$. Let
$\beta_1, \ldots , \beta_m$ be the roots of $f$ in $E$ and 
$L$ be the minimal $q$-polynomial of $f$ over $F$.
Let $V$ be the root space of $f$, consisting of all the roots of $L$.
 It is clear that $V$ is a module for $G$ over $\F_q$,
 of dimension $m(f)$.

 There is an abstract permutation module $W$ for $G$ 
with basis $w_1, \ldots ,w_m$, 
 where $G$ acts on the basis elements as it permutes the roots of $f$  in the Galois action.
 There is a $G$-module homomorphism 
 $\psi : W\longrightarrow V$ which maps each basis element to the corresponding root,
 i.e., $\psi (w_i)=\beta_i$.
 The kernel of $\psi$ is a $G$-submodule of $W$, and the image of 
 $\psi$ is $V$ by definition.
 Another way to say this is that $f$ gives rise to a short exact sequence 
 of $G$-modules
 \[
 0 \longrightarrow \ker (\psi) \longrightarrow W  \stackrel{\psi}{\longrightarrow} V
  \longrightarrow 0.
 \]

 The kernel of $\psi$ is determined by all the linear dependence relations among the roots of $f$,
 with coefficients in $\F_q$.
 The most obvious example occurs when the sum of the roots of $f$ is 0, and then the element
 $w_1+\cdots +w_m$ 
  is in the kernel of $\psi$.
  It may happen that $\ker (\psi)$ is 1-dimensional, with basis vector $w_1+\cdots +w_m$, 
  in which case $\dim (V)=m-1$, but if we can find 
  a polynomial $f$ with $\dim (V)<m-1$ then we have found a submodule (namely $\ker (\psi)$) of
  the permutation module which has dimension larger than 1.
  This gives another motivation for finding polynomials $f$ where $m(f)$ is significantly smaller
  than $\deg (f)$: they may give rise to interesting $G$-modules that appear as
  submodules of the permutation module.
  Other authors have considered linear relations among roots, mainly in characteristic 0, 
see \cite{BDEPS} or \cite{Gi} for example.

We close this section with some observations about $G$-modules and spaces of roots.
One may wonder which $G$-modules over $\F_q$ can occur as the space of roots
of a polynomial over $F$.

\begin{lemma}\label{not_hom_image}
Let $F$ contain $\F_q$. 
Any $G$-module which is not a homomorphic image of the 
$m$-dimensional permutation module $W$
cannot be the space of roots of any polynomial $f\in F[x]$ of degree $m$ that has Galois group $G$.
\end{lemma}

\begin{proof}
Any polynomial $f$ of degree $m$ gives rise to a module homomorphism
$\psi : W \longrightarrow V$ and a short exact sequence
 \[
 0 \longrightarrow \ker (\psi) \longrightarrow W  \stackrel{\psi}{\longrightarrow} V  \longrightarrow 0.
 \]
 Therefore the space of roots $V$ is a homomorphic image of $W$.
\end{proof}


As some sort of converse to this lemma we have the following lemma and corollary.

\begin{lemma}\label{anyirred}
Let $G$ be a Galois group over $F$. Then any  $\F_qG$-submodule of the regular
module $\F_qG$ is isomorphic to the space of roots of some $q$-polynomial over $F$ with Galois group 
a homomorphic image of $G$. 
\end{lemma}

\begin{proof}
Let $f\in \F[x]$ be irreducible and have Galois group $G$.
By the normal basis theorem there exists an irreducible polynomial $g$ of degree
$|G|$ with the same splitting field as $f$ and Galois group $G$. 
Let $L$ be the  minimal linearization of $g$.
Then $L$ has $q$-degree $|G|$, as the roots of $g$ are linearly independent over $\F_q$, 
and its space $S$, say, of roots is isomorphic as an $\F_qG$-module to the regular module $\F_qG$.

Let $U$ be an $\F_qG$-submodule of the regular module. Then by a slight abuse of notation, we may consider
$U$ to be a subspace of $S$.
Consider the polynomial
\[
L_U=\prod_{u\in U} (x-u),
\]
where we take the elements $u$ to be in the subspace of $S$ isomorphic to $U$ as
an $\F_qG$-module.
Then $L_U$ has coefficients in the splitting field of $f$ over $F$ and as the roots
of $L_U$ are the union of $G$-orbits, these coefficients are fixed elementwise by $G$
and thus lie in $F$. 

By Lemma \ref{lin}, $L_U$ is a $q$-polynomial in $F[x]$ and it divides $L$. Its space
of roots is isomorphic to $U$ as an $\F_qG$-module. Finally, the splitting field
of $L_U$ is contained in that of $L$, and thus that of $f$. Consequently,
the Galois group of $L_U$ is a homomorphic image of $G$. 
\end{proof}

\begin{cor} \label{irreducible_module}
Let $G$ be a Galois group over $F$. Then any  irreducible $\F_qG$-module is isomorphic  to the space of roots of some $q$-polynomial over $F$ with Galois group 
a homomorphic image of $G$. 
\end{cor}

\begin{proof}
This follows from Lemma \ref{anyirred}, since it is well known that any irreducible
$\F_qG$-module is isomorphic to a submodule of the regular module.
\end{proof}

\section{Possible Minimal Linearizations}\label{possible_min}

One might wonder which $q$-polynomials can occur as a minimal linearization, and if there is 
any hope of a classification.
We will next show that, when we take the field $F$ to be
$\F_q(t)$,  then under a mild hypothesis, which cannot be omitted,  any monic $q$-polynomial with coefficients
in $\F_q[t]$ is
the \emph{minimal} linearization of at least one of its irreducible factors. 
This shows that minimal linearizations can be arbitrarily complicated, and there is little hope of
a characterization. 

Before the  result we state a lemma from \cite{BSEM} (Proposition 2.11) that we will need.

\begin{lemma} \label{cyclic_element}
Let
\[
L=x^{q^n}+a_{n-1}(t)x^{q^{n-1}}+\cdots +a_1(t)x^q+a_0(t)x,
\]
be a monic $q$-polynomial of $q$-degree $n$ whose coefficients $a_i(t)$ are elements
of the polynomial ring $\F_q[t]$. Let $G$ be the Galois group of $L$ over $\F_q(t)$.
Suppose that $a_0(\lambda)\neq 0$ for some element $\lambda$ of $\F_q$.  
Then there is an element in $G$ that acts as a cyclic linear
transformation on the space $V$ of roots of $L$ whose minimal polynomial for this action is 
\[
x^n+a_{n-1}(\lambda)x^{n-1}+\cdots + a_1(\lambda)x+a_0(\lambda).
\]
\end{lemma}

To be precise, the statement in  \cite{BSEM} does not mention the minimal
polynomial, it concerns the characteristic polynomial. However one can 
see that these are equal under the stated hypotheses. A key observation, which 
underlies the principle and is quite easy
to prove, is that if we have a $q$-polynomial 
\[
P=x^{q^n}+b_{n-1}x^{q^{n-1}}+\cdots +b_1x^q+b_0x
\]
in $\F_q[x]$, then its associate
\[
x^n+b_{n-1}x^{n-1}+\cdots +b_1x+b_0
\]
is the minimal polynomial of the Frobenius $q$-power map acting on the space of roots of $P$.


The following theorem 
ensures that almost every monic $q$-polynomial over $\F_q[t]$ 
occurs as a minimal $q$-polynomial of an irreducible polynomial. 

\begin{thm}\label{anyLpoly}
Let $L=\sum_{i=0}^n a_i(t) x^{q^i}$ be a monic $q$-polynomial  with coefficients in  $\F_q[t]$  such that $t^q-t$ does not divide $a_0(t)$.
Then $L$ has an irreducible factor $f$ such that $L$ is the minimal $q$-polynomial
of $f$.
\end{thm}

\begin{proof}
Let $L$ be a monic $q$-polynomial  with coefficients in  $\F_q[t]$  such that $t^q-t$ does not divide $a_0(t)$.
Let $E$ be the splitting field of $L$ over $\F_q(t)$.
Let $G$ be the Galois group of the extension $E:\F_q(t)$.  
Let  $V \subseteq E$ be the space of roots of $L$. Let $m=\dim V$ = $q$-degree of $L$.
By Lemma \ref{cyclic_element} there is a $\sigma\in G$ that acts 
as a cyclic linear
transformation on  $V$.
This means that there exists $\alpha\in V$ such that the elements
$\sigma^i (\alpha)$ for $i=0,1,\ldots, m-1$ are a basis for $V$.

 The minimal polynomial of $\alpha$ over $\F_q(t)$
will have the $\sigma^i (\alpha)$  among its roots.
This polynomial 
 is irreducible over $\F_q(t)$, and has minimal linearization $L$. (We remark
 that, more precisely, the minimal polynomial of $\alpha$ has coefficients in $\F_q[t]$, since this
 is true of $L$.)
 \end{proof}
 
 Now we want to explain our remark that the hypothesis concerning $t^q-t$ cannot in general be omitted
 in Theorem \ref{anyLpoly}. We do this by providing an example for $q=3$.
 
 \begin{thm} \label{counterexample} There exists a $3$-polynomial $L(x)$ of $3$-degree $4$ in $\F_3(t)[x]$ that is not
 the minimal $3$-polynomial  of any of its irreducible factors.
 The coefficient of $x$ in $L(x)$ is divisible by $t^3-t$.
 \end{thm}
 
 \begin{proof}
 We outline the details, which involve some  calculations. We start
 with the polynomial $f=(x^2+x+t)(x^3+x+t)$ in $\F_3(t)[x]$, which is the product
 of two irreducible polynomials. Using Magma, we calculate
 the minimal 3-polynomial, $L$, say of $f$, which has 3-degree 4. We will not write down $L$
 explicitly, as its coefficients are quite large polynomials in $t$.
 
 Again using Magma, we find the irreducible factors of $L$ over $\F_3(t)$. 
 There are three  different irreducible factors of degree 1 (which are
 $t$, $t-1$, $t-2$), six irreducible factors of degree 2,
 six irreducible factors of degree 3, three irreducible factors of degree 4 and six irreducible 
 factors of degree 6.
 
 We have to show that $L$ is not the minimal 3-polynomial of any of these irreducible
 factors. This we do by eliminating each in turn using Magma. By the upper bound
 in Lemma \ref{bounds}, we need only
 consider the irreducible factors of degree 4 or 6. We find that the minimal
 $3$-polynomial of each of those factors has 3-degree 3, and hence is never $L$, as required.
 \end{proof}

While the  occurrence of the polynomial $t^q-t$ as a divisor of the $x$-coefficient of
the minimal $q$-polynomial of a polynomial in $\F_q[t][x]$ is presumably hard to predict, we will describe one case where it does happen.

Let $E$ be a Galois extension of $\F_q(t)$ and $G=Gal(E:\F_q(t))$.  
Let $\alpha_1, \ldots , \alpha_{|G|}$
be a normal basis for $E$ over $\F_q(t)$. By definition, this basis consists of all
the $G$-conjugates $\sigma_i(\alpha)$ of some element $\alpha$ in $E$, where the $\sigma_i$
run over $G$. Certainly, $\alpha$ is algebraic over 
$\F_q(t)$ but it may not be algebraic over $\F_q[t]$ (meaning that it is not necessarily an algebraic integer 
in $E$). However, it is an elementary result from the theory of Dedekind domains
that we can find a polynomial $g(t)$, say, in $\F_q[t]$ such that $\beta=g(t)\alpha$ 
is an algebraic integer. Setting $\beta_i=\sigma_i(\beta)$, for all $i$, 
the $\beta_i$ are all algebraic integers in $E$. The $\beta_i$ are also linearly independent over $\F_q(t)$, and hence over $\F_q$. 
Thus their $\F_q$-span, $V$, say, 
has dimension $|G|$. 

Let $P$ be the monic $q$-polynomial over $\F_q(t)$
whose roots are the elements of $V$. This has coefficients in $\F_q[t]$, as the $\beta_i$
are algebraic integers in $E$. 
Let   $f$ be the minimal polynomial of $\beta=\beta_1$ over $\F_q(t)$. Likewise,
the coefficients of $f$ are elements of $\F_q[t]$. 
Then $P$ is the minimal $q$-polynomial of $f$, and $m(f)=\deg f=|G|$, because the
roots of $f$ are linearly independent over $\F_q$.

Suppose $t^q-t$ does not divide the coefficient of $x$ in $P$.
Then by Lemma  \ref{cyclic_element} there is a Frobenius element $\sigma$ in 
$G$ that acts on $V$ as a cyclic element and hence its minimal polynomial  $\pi$ has degree $|G|$.
However, if $\sigma$ has order $r$ 
then $x^r-1$ must be divisible by $\pi$. Hence $|G|\le r$. As $r$ divides $|G|$, this implies
$r=|G|$ and $G$ is cyclic and generated by $\sigma$.

We can summarize this argument as follows. 

\begin{thm} \label{noncyclic_case}
Let $E$ be a Galois extension of $\F_q(t)$ whose Galois group $G$ is not cyclic. Let $\alpha$ be any
element of the ring of integers of $E$ whose $G$-conjugates form a normal basis
of $E$ over $\F_q(t)$ (such elements exist, as just outlined). Let $f$ be the minimal polynomial
of $\alpha$ over $\F_q(t)$ and let $P$ be its minimal $q$-polynomial. Then both $f$
and $P$ are monic polynomials with coefficients in $\F_q[t]$, of degree $|G|$ and $q^{|G|}$, respectively, and $t^q-t$ divides the $x$-coefficient of $P$.
\end{thm}

\section{Permutation modules related to $S_m$}\label{sect_sm}

\noindent If we picked a random polynomial $f$, what would $m(f)$ be? Here we discuss this question, but we make no attempt to quantify our observations.

Let $W$ be the permutation module of dimension $m$ over $\F_q$ for the symmetric group $S_m$.
We take $W$ to have a basis $w_1$, \dots, $w_m$ permuted according to the action of $S_m$
on the $m$ indices. Let
\[
e=w_1+\cdots + w_m.
\]
This element is clearly fixed by $S_m$ and, up to scalar multiples, it is the unique
nonzero element in $W$ fixed by $S_m$. This holds for any transitive subgroup of degree $m$.
We define $W_1$ to be $\{ \lambda e : \lambda \in \F_q \}$,
 a subspace of $W$ which is $S_m$-invariant and has dimension one in $W$.

Let $W_0$ be the subspace of $W$ consisting of all sums of the form
\[
\lambda_1 w_1+\cdots +\lambda_m w_m,
\]
where the $\lambda_i$ are elements of $\F_q$ and
\[
\lambda_1+\cdots+\lambda_m=0.
\]
It is easy to see that $W_0$ is $S_m$-invariant and has codimension one in $W$. We call
$W_0$ the deleted permutation module for $S_m$.

The following is well known, but we include a proof for completeness.

\begin{lemma} \label{submodules_theorem}
Let $U$ be a proper $S_m$-invariant $\F_q$-subspace of the permutation module $W$ for $S_m$.
Then $U=W_1$ or $U=W_0$.
\end{lemma}

\begin{proof}
Suppose that $U\neq W_1$. Then $U$ contains an element
\[
u= \lambda_1 w_1+\cdots +\lambda_m w_m,
\]
where not all the $\lambda_i$ are equal, say $\lambda_j\neq \lambda_k$ for some $j$ and $k$
with $1\leq j<k\leq m$.

Let $\tau$ be the transposition $(j\ k)$ in $S_m$. Then we find that 
\[
\tau(u)-u=(\lambda_k-\lambda_j)(w_j-w_k)
\]
is in $U$. Thus since $\lambda_j\neq \lambda_k$, the element $w_j-w_k$ is in $U$.
Then, given the $m$-fold transitivity of $S_m$, we find that $w_a-w_b$ is in $U$
for all indices $a$ and $b$ with $a\neq b$. Since these elements clearly span
$W_0$, we see that $U=W_0$.
\end{proof}

This lemma implies that there are no interesting submodules of the permutation module
for $S_m$. 
 In later sections we will see that there are some interesting submodules
for some other highly transitive groups.

We consider next a monic polynomial $f$ of degree $m$ in $F[x]$. We call the coefficient
of $x^{m-1}$ in $f$ the \emph{trace} of $f$. The vanishing, or otherwise, of the trace
of $f$ is relevant to the study of $m(f)$, for the trace of $f$ is zero
if and only if the sum of the roots of $f$ in some splitting field over $F$ is zero.

The following is a reasonably generic special case, 
because a randomly chosen polynomial has Galois group $S_m$
(see \cite{EntinPopov} for example, for the case $F=\F_q(t)$).
This  allows us to produce examples
of polynomials $f$ with $m(f)=\deg(f)$, meeting the bound of Lemma \ref{bounds}.

\begin{thm} \label{the_value_of_m(f)}
Let $f$ be a monic irreducible polynomial of degree $m$ over $F$. Suppose that the Galois
group over $F$ is isomorphic to the symmetric group $S_m$. Then
$m(f)=m-1$ is $f$ has trace zero, and $m(f)=m$ otherwise.
\end{thm}

\begin{proof}
Let $V$ be the $\F_q$-subspace spanned by the roots $\beta_i$, say, of $f$ in some splitting field over $F$. Let $W$ be the abstract permutation module for $S_m$ over $\F_q$ and let
$\psi$ be the evaluation map from $W$ to $V$, under which a basis element $w_i$ is sent
to the root $\beta_i$ for all $i$. 
The kernel $\ker \psi$ of $\psi$ is an $S_m$-invariant subspace of $W$, and it is certainly
not the whole of $W$. 

Suppose that $f$ has trace zero. Then $\ker \psi$ contains $W_1$, since the sum of the roots
is zero. We must then show that $\ker \psi=W_1$. Suppose that $\ker \psi$ is not
$W_1$. It follows from Lemma \ref{submodules_theorem} that $\ker \psi=W_0$. But this implies
that all the roots of $f$ are equal, which is not the case. We deduce that
$\ker \psi$ is one-dimensional and thus $V=\psi(W)$ has dimension
$m-1$. This of course implies that $m(f)=m-1$.

Now suppose that $f$ has nonzero trace. It follows that $W_1$ is not contained in
$\ker \psi$. Then Lemma \ref{submodules_theorem} implies that 
either $\psi$ is a monomorphism or $\ker \psi=W_0$. The second possibility is excluded,
as it implies that the roots of $f$ are all equal. Therefore, $\psi$ is a monomorphism
and $V=\psi(W)$ has dimension $m$. Hence $m(f)=m$ in the nonzero trace case.
\end{proof}

We conclude from this theorem that a randomly chosen polynomial $f$
is highly likely to have $m(f)=\deg(f)$.

\section{Golay code and Steiner system on 24 points}\label{golay}

\noindent We present a construction of the $\F_2$ Golay code of length 24, and the Steiner system on 24 points. This construction makes several aspects of these structures very visible.
One interesting aspect of this construction is that the two objects are constructed together.

We start with the polynomial $f=x^{24}+x+t \in \F_2(t)[x]$, which is irreducible over $\F_2(t)$.
We do not assume anything about the Galois group of $f$, and in particular we do not assume
that the Galois group is $M_{24}$.
We do assume that  the minimal 2-polynomial $L$ of $f$ is
the polynomial given in Table 1,
and that the irreducible factors of $L$ over $\F_2(t)$ have degrees 
$1, 24, 276, 1771, 2024$.
These facts are quickly established by a computer algebra package such as Magma, and
can in principle be done by hand (e.g. $L$ is calculated by hand in \cite{CMT}).
With these assumptions, we shall give a simple construction of the Golay code and Steiner system.

Let $W$ be the permutation module of dimension 24, and let $V$ be the  space of roots of $f$.
The roots of $L$ are the elements of $V$.
There is a surjective module homomorphism $\psi: W\longrightarrow V$ with a 
12-dimensional kernel.

We shall next show that the Galois group of $f$ is 5-transitive.
This can be proved by elementary methods, as we now show for completeness.
This result is stated in \cite{CMT}  without proof, and is also stated in  \cite{AY} 
where the proof uses
the classification of finite simple groups.

\begin{lemma}
The Galois group of $f$ over $\F_2(t)$ is 5-transitive.
\end{lemma}

\begin{proof}
Let $G$ be the Galois group of $f=x^{24}+x+t$ over $\F_2(t)$.
We use the specialization theorem of Dedekind.
Specializing $t$ to the elements $a$ of $\F_8$, and factoring $x^{24}+x+a$ into irreducible factors, 
shows that the Galois group $G$ contains elements with the following cycle structure:
$(23)(1)$,
$(7)^3 (1)^3$,
$(11)^2 (1)^2$.

Specializing $t$ to the elements $a$ of $\F_{64}$, and factoring $x^{24}+x+a$ into irreducible factors, 
shows that the Galois group contains elements with the following cycle structure:

$(12)^2$, whose cube is $(4)^6$, and sixth power is $(2)^{12}$, and

$(1)^2 (2) (4) (8)^2$, whose fourth power is $(1)^8 (2)^8$.

The elementary argument in the proof of Theorem 2, section 3, chapter 20 of \cite{MS} shows that 
these cycle types imply that $G$ is 5-transitive.
\end{proof}

We remark that it is relatively easy to prove that the group $G$ is 4-fold transitive using cycle types but it requires more ingenuity to prove that is is 5-fold transitive.

We order the roots of $f$ as $\beta_1, \ldots , \beta_{24}$
 and  we fix a corresponding  basis $w_1, \ldots ,w_{24}$ of $W$.
 There is a surjective module homomorphism $\psi: W\longrightarrow V$ defined by 
 $\psi (w_i)=\beta_i$.

\begin{lemma}\label{evenweight}
An odd number of roots of $f$ cannot sum to 0.
\end{lemma}

\begin{proof}
Continue the above notation. 
The module homomorphism $\psi: W\longrightarrow V$ has a 
12-dimensional kernel, which we denote by $C$.

There is a $G$-invariant  linear functional $\lambda : W \longrightarrow \F_2$
defined by 
\[
\lambda (\sum_{i=1}^{24} a_i w_i ) = \sum_{i=1}^{24} a_i.
\]
The kernel of $\lambda$ is the deleted permutation module.
We want to prove that $\lambda$ is identically 0 when restricted to $C$.
This will prove the lemma.

Suppose not. Then $\lambda$ has an 11-dimensional kernel on $C$.
As shown in the previous proof, $G$ has an element of order 23, call it $\sigma$.
Then $\sigma$ acts on $\ker (\lambda)$, and in fact $\sigma$ must act
irreducibly on $\ker (\lambda)|_C$. However, the vector 
$e=\sum_{i=1}^{24}  w_i$ is in $\ker (\lambda)|_C$, and 
generates a 1-dimensional submodule for $\sigma$, a contradiction.
\end{proof}

\begin{lemma}\label{sumfour}
Four roots of $f$ cannot sum to 0.
\end{lemma}

\begin{proof}
Suppose $\beta_i+\beta_j+\beta_k+\beta_\ell=0$.
By the 4-transitivity of $G$ we could choose an element of $G$ that fixes 
$\beta_i$ and $\beta_j$ and $\beta_k$, and doesn't fix $\beta_\ell$.
Say it maps $\beta_\ell$ to $\beta_m$.
Applying this element gives an equation 
$\beta_i+\beta_j+\beta_k+\beta_m=0$ where $m\not= \ell$.
This implies $\beta_\ell = \beta_m$, a contradiction.
\end{proof}

\begin{cor}
Sums of two roots of $f$ are all distinct.
\end{cor}

Since $\binom{24}{2}=276$ the polynomial whose roots are all sums of two roots of $f$
must have degree 276. This accounts for the irreducible factor of $L$ of degree 276.

\begin{lemma}\label{sumsix}
Six roots of $f$ cannot sum to $0$.
\end{lemma}

\begin{proof}
Suppose $\beta_i+\beta_j+\beta_k+\beta_\ell+\beta_m+\beta_n=0$
where the roots involved are all different.
Choose a root $\beta_r$ which is not equal to any of these six roots.
By the 5-transitivity of $G$ we can choose an element $\sigma\in G$ that fixes 
$\beta_i$ and $\beta_j$ and $\beta_k$ and $\beta_\ell$, 
and has $\sigma (\beta_m)=\beta_{r}$.
Let $\beta_s$ be $\sigma(\beta_n)$.
Applying $\sigma$ and adding the two equations gives 
$\beta_m+\beta_n=\beta_{r}+\beta_{s}$.
Since $r\not=m$, $r\not=n$, $r\not= s$,  we have a contradiction to Lemma \ref{sumfour},
or to the distinctness of the roots.
\end{proof}

\begin{cor} \label{sums_of_three_roots}
Sums of three roots of $f$ are all distinct.
\end{cor}

Since $\binom{24}{3}=2024$ the polynomial whose roots are all sums of three roots of $f$
must have degree 2024. This accounts for the irreducible factor of $L$ of degree 2024.

Since Corollary \ref{sums_of_three_roots} is important in our exposition, we propose to give
an alternative proof which only uses the 3-transitivity of $G$ and not the full force of
5-transitivity.

\begin{lemma} \label{degree_2024_factor}
The roots of the degree $2024$ factor of $L$ are the sums of three different roots of $f$.
Thus sums of three different roots of $f$ are all distinct and there are $2024$ such sums.
\end{lemma}

\begin{proof}
Let $\beta$ be the sum of three different roots of $f$. We know 
 that $\beta$ is not equal to a root of $f$ (by Lemma \ref{sumfour}) nor equal to a sum of two
roots of $f$ (by Lemma \ref{evenweight}). Since $\beta$ is a root of $L$, 
it follows from the factorization of $L$ into irreducibles that
$\beta$ is either a root of the irreducible factor of degree 1771 or that of degree 2024.

If $\beta$ is a root of the polynomial of degree 2024, then there are 2024 different sums of three
roots, since they are all in a single Galois group orbit by irreducibility, and that is what we want.

So let us assume for the sake of contradiction
that $\beta$ is a root of the degree 1771 polynomial. The $G$-orbit of $\beta$ has size 1771 as the polynomial is irreducible.
Let $H$ be the stabilizer subgroup of $\beta$ in $G$. We have
$|G|=1771\times |H|$.

As we know that $G$ acts 3-transitively on the roots of $f$, it acts transitively on the unordered subsets
of the three roots, of which there are
\[
24 \times 23 \times 22/6 =4 \times 23 \times 22
\]
members.
Let $M$ be the subgroup of $G$ fixing such an unordered 3-subset.
Then $|G|=4 \times 23 \times 22 \times |M|$, by transitivity.

Now by the transitive action on 3-subsets of roots of $f$, 
we may assume that $\beta$ is the sum of the three roots in this fixed 3-subset. It follows that
$M$ fixes $\beta$. 
This implies that $M$ is a subgroup of $H$ and thus $|M|$ divides $|H|$.

But we have 
\[
4 \times 23 \times 22\times |M|=1771\times |H|
\]
 and hence we obtain an integer
\[
|H|/|M|=8 \times 11 \times 23/1771,
\]
which is impossible since $1771=7 \times 11 \times 23$.
This proves the required result. 
\end{proof}

So far we have accounted for the roots of $f$, sums of two roots, and sums of three roots, 
which correspond to the factors of $L$ having degrees 24, 276, 2024, respectively.
Apart from $x$ there is one irreducible factor of $L$ remaining, of degree 1771.
This factor must correspond to sums of four roots, because a sum of four roots of $f$ cannot
equal a root of $f$ (by Lemma \ref{evenweight}) or a sum of two roots of $f$ 
(by Lemma \ref{sumsix}) or a sum of three roots of $f$  (by Lemma \ref{evenweight}).

We remark that Todd, in Section 2.1 of \cite{T2}, reaches similar conclusions about the action of
$M_{24}$ on the 12-dimensional Todd module over $\F_2$, but he does not work
in the context of roots of polynomials.

Next, we deal with sums of eight roots.

\begin{lemma}
The roots of $f$ form a Steiner system $S(5,8,24)$.
There are $759$ sets of eight roots that sum to $0$.
\end{lemma}

\begin{proof}
Consider a sum of five roots of $f$, call it $S$.
This $S$ cannot be 0 by Lemma \ref{evenweight}.
However $S$ must be a root of one of the irreducible factors of $L$,
which have degrees 1, 24, 276, 1771, and 2024.
But $S$ cannot equal a root of $f$ by Lemma \ref{sumsix}, and cannot 
equal a sum of two roots or a sum of four roots by Lemma \ref{evenweight}.
The only option is that $S$ is equal to a sum of three roots, and $S$
is a root of the degree 2024 factor.

Since $S$ is equal to a sum of three roots, we obtain a sum of eight roots which is 0.
Furthermore, $S$ can only equal a sum of three roots in one way, because if 
$S=T_1$ and $S=T_2$ (where $T_1$ and $T_2$ are sums of three roots)
then $T_1=T_2$, and then $T_1+T_2=0$ which
contradicts Lemma \ref{sumsix}.
Therefore, any set of five roots determines a unique set of three further roots such that
the sum of all eight roots is 0.

The number of sets of eight roots that sum to 0 is therefore
$\binom{24}{5}/\binom{8}{5}=759$.
\end{proof}

A set of eight roots of $f$ that sum to 0 is called an octad.
This proof has shown that the roots of $f$ form a Steiner system, 
where the blocks are the octads.
This is a realization of the Steiner system $S(5,8,24)$, together with an action of
the Galois group of $f$.

\begin{lemma}
The kernel of $\psi :  W \longrightarrow V$ is isomorphic as a vector subspace
of $(\F_2)^{24}$ to the binary Golay code of length $24$.
\end{lemma}

\begin{proof}
We consider binary vectors in $(\F_2)^{24}$ and we label the coordinates by 
the basis $w_1, \ldots ,w_{24}$, or equivalently by 
the roots of $f$. We construct a 12-dimensional subspace of $(\F_2)^{24}$
from the kernel of $\psi$ by putting
a 1 in position $w_i$ for each element $\Sigma_i w_i$ of the kernel.
Previous lemmas show that none of these vectors has weight less than 8, and that
759 of these vectors have weight 8.
We have constructed a binary $[24,12,8]$ code with 759 codewords of weight 8.
By various standard arguments, the code must be the Golay code.
\end{proof}

The construction of the Golay code and Steiner system is complete.
Now we shall return to sums of four roots. The remainder of this section is not necessary
for the construction, however we wish to complete the picture.

A sum of four roots cannot equal a root or a sum of three roots by Lemma \ref{evenweight},
and cannot equal a sum of two roots by Lemma \ref{sumsix}.
Therefore any sum of four roots must be a root of the degree 1771 factor of $L$.
It follows by Galois action that all   roots of the degree 1771 factor must be 
sums of four roots.

Let $A$ be a subset of four roots of $f$. We want to describe a way of finding octads
containing $A$. Let $\beta$ be any of the twenty roots not in $A$. Then $A \cup \{\beta\}$ is a subset of five roots and hence it lies in a unique octad, $D$, say, by the Steiner property.
Then $D=A\cup A'$ for some subset $A'$ of four roots with no root in common with
$A$. If we had replaced $\beta$ by any root $\gamma$ in $A'$, this procedure would lead
to the same octad $D$. Thus running through the twenty possibilities for $\beta$,
we obtain five different octads, $D_1$, \dots, $D_5$, where $D_i\cap D_j=A$, for $i\neq j$.

Write $D_i=A\cup A_i$, for $1\leq i\leq 5$. The sum of the roots in $A_i$ 
is equal to the sum of the roots in $A$ for all $i$.
We now have five additional ways of expressing a sum of four roots as a sum of four roots, and so
with the original sum of roots in $A$, we get six ways.

Let $X_4$ be the set of all sums of four distinct roots of $f$, i.e., the set of roots of the
degree 1771 irreducible factor of $L$.
This argument shows that any element of $X_4$ can be expressed as a sum of 
four roots of $f$ in at least 6 different ways. 
There are 
\[
\binom{24}{4} = 10626=6\cdot 1771
\]
possible ways of choosing four distinct roots.
Since there are at least 6 different ways of expressing an element of $X_4$
as a sum of four distinct roots, there are exactly 6 ways.

We can see this in terms of properties of the Galois group $G$. Let $\alpha_i$, $1\leq i\leq 4$,
be the roots in our original subset $A$. Then $\alpha_1+\cdots +\alpha_4$ is in $X_4$. Let $\Gamma$ be the subgroup of $G$ that fixes
$\alpha_1+\cdots +\alpha_4$ and let $H$ be the subgroup of $G$ that maps $A$ into itself.
Then $H$ is a subgroup of $\Gamma$ and we find that $|\Gamma:H|=6$. $\Gamma$ 
transitively permutes the six
subsets $A$, $A_1$, \dots, $A_5$ among themselves and $H$ is the stabilizer of a subset
in this action.

A set of  six subsets of four roots with the same sum is called a sextet.

Finally, we show that $G=Gal(f)$ is isomorphic to $M_{24}$.
The argument about the order of $G$ in 
 the proof of Theorem 4, section 4, chapter 20 of \cite{MS} goes through. 
 The proof uses the lemmas and discussion above.
 This implies that  $G$ is a 5-transitive group of order 244823040.
 This is the order of $M_{24}$.
It is not hard to then prove that $G$ is isomorphic to $M_{24}$.

We have constructed a 12-dimensional submodule $C$ of the 24-dimensional permutation
module $W$ and we have proved that this submodule is the Golay code. 
This submodule is known as the Golay module for $M_{24}$.
 The quotient $W/C$ is the space of roots, also a 12-dimensional module - see next section.

We remark that a result known as Todd's Lemma becomes trivial in our context. 

\begin{lemma}[Todd's Lemma]
Suppose $D$ and $D'$ are octads that have an intersection of cardinality $4$.
Then the symmetric difference $D\ \Delta\ D'$ is an octad.
\end{lemma}

\section{Permutation modules related to $M_{24}$}\label{m24}

\noindent Let $W$ be the permutation module of dimension $24$ over $\F_2$ for the 
Mathieu group $M_{24}$.
We take $W$ to have a basis $w_1$, \dots, $w_{24}$ permuted according to the action of $M_{24}$
on the $24$ indices. Let
\[
e=\sum_{i=1}^{24}  w_{i}.
\]
This element $e$ is clearly fixed by $M_{24}$ and generates a 1-dimensional
submodule of $W$, which we will denote by $W_1$.

Lemma 1.2.1 of \cite{I} states that $W$ is a uniserial module with unique composition series
\[
\{0\} < W_1 < C < W_{0} < W,
\]
where $W_{0}$ denotes the 23-dimensional submodule of $W$ consisting of all
sums of an even number of basis vectors and $C=\ker \psi$, as described in the previous section.
The proof of this uniserial property is not obvious, however we will use the result.
Clearly then $C$ has a 1-dimensional submodule, with an 11-dimensional quotient.
By Lemma \ref{not_hom_image} 
the Golay module $C$ cannot be realized as the module of roots of a degree 24 polynomial
with Galois group $M_{24}$, 
because it is not a quotient of  $W$ (invoking the uniseriality of $W$).
As $M_{24}$ is simple, the Golay module  can however be realized as the module of roots of a higher degree polynomial,
as Lemma \ref{anyirred} shows.

The quotient $W/C$ is the known as the Todd module \cite{T1}
and is isomorphic to the dual of the Golay module. Therefore we have realized the Todd module
as the space of roots of $f$.
The quotient $W_{0}/C$ is called the irreducible 11-dimensional Todd module. Likewise,
$C/W_1$ is called the irreducible 11-dimensional Golay module and is the dual of the
irreducible Todd module.

In summary, for $M_{24}$, we can't get any module of dimension 11 as the space of roots for any 
linearization of a degree 24 polynomial with Galois group  $M_{24}$.
We can get a module of dimension 12 (the Todd module) but not its dual, the Golay module.  
We could possibly find a $g$ that has $m(g)=11$, whose minimal
linearization has roots which are the elements of this submodule.
Indeed, the polynomial $f_1$ in Table \ref{tab:experiments1}
is such a polynomial, it has degree 276.
There is no simple explanation for the existence of the Golay module as a 
submodule of the permutation module, nor of the properties of the quotient, the Todd module.

The degrees of the irreducible factors of $L$ imply that $M_{24}$
acts on the space of roots with orbits of sizes 1, 24, 276, 1771, 2024.
The group $M_{24}$ also acts on the Golay module, however the orbits have sizes
1, 1, 759, 759, and 2576.
This is a situation where the orbit sizes on a module and its dual are different, but it is always
the case that a module and its dual have the same number of orbits.

One may ask if a polynomial $f$ of degree 24 with Galois group $M_{24}$ is forced to have $m(f)=12$. This is not the case, as we now intend to demonstrate. In order to achieve this, we need to assemble
some facts from representation theory. These are well known, but we provide proofs for want
of an obvious reference.

Let $G$ be a finite group and let $M$ be an $\F_qG$-module. We say that $M$ is a cyclic
$\F_qG$-module if there is an element $\alpha$ in $M$ such that the $\F_q$-span of the elements
$\sigma(\alpha)$, $\sigma\in G$, equals $M$. Such an element $\alpha$ is called a generator of $M$.
It is straightforward to prove that the cyclic $\F_qG$-modules are precisely
the homomorphic images of the regular module $\F_qG$.

Given an $\F_qG$-module $M$, there is the associated dual $\F_qG$-module $M^*$, which
consists of the space of $\F_q$-linear maps from $M$ into $\F_q$. Given an $\F_qG$-submodule, $U$, say, of
$M$, we let $A(U)$ denote the annihilator of $U$ in $M^*$. It consists of all linear
maps $\phi$ in $M^*$ that satisfy $\phi(u)=0$ for all $u\in U$. $A(U)$ is an $\F_qG$-submodule
of $M^*$ and we have
\[
\dim U+\dim A(U)=\dim M.
\]
Clearly, if $U<W$ for some $\F_qG$-submodule $W$ of $M$, $A(W)<A(U)$.

\begin{lemma} \label{unique_max_and_min}
Let $M$ be an $\F_qG$-module that contains both a unique minimal and a unique maximal
submodule. Then $M^*$ contains both a unique minimal and a unique maximal
submodule. Moreover, both $M$ and $M^*$ are cyclic $\F_qG$-modules.
\end{lemma}

\begin{proof}
Let $M_0$ and $M_1$ be the unique minimal and unique maximal submodules of $M$, respectively.
Let $\alpha$ be any element of $M$ not contained in $M_1$. Then it is easy to see that
$\alpha$ is a generator for $M$ as an $\F_qG$-module and hence
$M$ is cyclic.

Furthermore, it is also easy to see that $A(M_1)$ is the unique minimal $\F_qG$-submodule 
of $M^*$ and $A(M_0)$ is the unique maximal submodule. Thus $M^*$ is also a cyclic
module.
\end{proof}

We remark that it is not clear that if we assume only that $M$ is a cyclic $\F_qG$-module,
then $M^*$ is also cyclic.

\begin{lemma} \label{submodule_of_regular}
Let $M$ be an $\F_qG$-module that contains both a unique minimal and a unique maximal
submodule. Then $M$ is isomorphic to an $\F_qG$-submodule of the regular $\F_qG$ module.
\end{lemma}

\begin{proof}
Lemma \ref{unique_max_and_min} implies that $M^*$ is cyclic. Hence there is an $\F_qG$ submodule,
$I$, say, of $\F_qG$ such that $F_qG/I\cong M^*$.

Now it is well known that there is a nondegenerate symmetric bilinear form, $B$, say, defined
on $\F_qG\times \F_qG$ that is $G$-invariant. Let $I^\perp$ denote the perpendicular subspace
of $I$ in $\F_qG$ with respect to $B$. Then general duality theory shows that
\[
I^\perp\cong (\F_qG/I)^*\cong (M^*)^*.
\]
But it is well known that $(M^*)^*$ is isomorphic to $M$ as an $\F_qG$-module, and thus 
$I^\perp$ provides the required submodule.
\end{proof}

We return to the permutation module $W$ for $M_{24}$ over $\F_2$ and recall that $W$ is uniserial.
Let $W^1$ be the quotient $W/W_1$. Since $W$ is uniserial, so also is $W^1$. It follows that
$W^1$ has a unique minimal and unique maximal submodule. Lemma \ref{submodule_of_regular}
thus implies that $W^1$ is isomorphic to a submodule of the regular module for $M_{24}$
over $\F_2$

\begin{thm}
Let $F$ be a field of characteristic $2$ having a Galois extension with Galois group $M_{24}$.
There exists an irreducible $f\in F[x]$ of degree $24$ having Galois group $M_{24}$ and $m(f)=23$.
\end{thm}

\begin{proof}
Let $E$ be a Galois extension of
$F$ with Galois group $M_{24}$. 
Then by the normal basis theorem, there is a vector subspace $V$ 
of dimension $|M_{24}|$ over $\F_2$ inside $E$ which is invariant under
$M_{24}$. Note that $V$ is isomorphic to the regular module for $M_{24}$ over $\F_2$.

Let $w_1$ be an element of the permutation basis for $W$ and $\alpha = w_1 +W_1$ in the quotient $W^1$. Now the orbit of $\alpha$ 
under the action of $M_{24}$ has size at most 24, since the orbit of $w_1$ under $M_{24}$ has size exactly 24. The stabilizer of $\alpha$ under the action of $M_{24}$ clearly contains
the stabilizer of $w_1$, which is the maximal subgroup $M_{23}$ of $M_{24}$.
Hence by maximality, $\alpha$ has exactly 24 conjugates in $W^1$ under $M_{24}$ action.

By Lemma \ref{submodule_of_regular}, there is a submodule of $V$ isomorphic to $W^1$.
Let $L$ be the monic 2-polynomial  of 2-degree 23 whose roots are the elements of this submodule.
 The Galois group of $L$ is
$M_{24}$ over F. There is then a root $\beta$ of $L$ whose orbit under $M_{24}$
has size 24 (we obtain $\beta$ from the element $\alpha$ just described).

Let $f$ be the minimal polynomial of $\beta$ over $F$. This has degree 24
because $\beta$ has 24 conjugates under the Galois group action.
So we have a polynomial $f$ of degree 24 over $F$ which divides $L$, and we want
to show that $L$ is the minimal 2-polynomial for $f$.

Let $M$ be the minimal 2-polynomial for $f$ over $F$. Then $M$ divides $L$ and so its space
of roots is a subspace of the space of roots of $L$, invariant
under the action of $M_{24}$ (in other words, a submodule of the space
of roots of $L$).

The space of roots of $M$ is the $\F_2$ span of the roots of $f$. Now $M_{24}$
contains a unique conjugacy class of subgroups of index 24 (this means
that there is essentially only one way in which $M_{24}$ acts on 24 points).

This means that as a module for $M_{24}$, the space of roots of $f$ is isomorphic
to a quotient of $W$. By Lemma 1.2.1 of \cite{I}, the nontrivial quotients of $W$ are the following:
$W^1$; a 12-dimensional module, 
$W/C$, which is the Todd module; and a 1-dimensional module $W/W_{0}$.
The space of roots of $M$ cannot be
one-dimensional. So, if $M$ is not $L$, its space of roots must be 12-dimensional
and isomorphic to the Todd module.

However, $W^1$ does not contain a 12-dimensional submodule
because $W^1$ is also uniserial (it is a quotient
of $W$), and the only submodules of $W^1$ have dimensions 11 and 22.
\end{proof}

This result means that the Golay code does not inevitably appear in this process.
We could have been unlucky in the choice of $f=x^{24}+x+t$  and have a linearization of large degree
that would not help us.
It remains an open problem to find explicitly a degree 24 polynomial $g$ with 
Galois group $M_{24}$ and $m(g)=23$.

We can also find degree 24 polynomials $g$ with Galois group $M_{24}$ and $m(g)=24$.
This is more straightforward. For example, working over $\F_2(t)$, both
$x^{24}+x+t$ and $x^{24}+x^{23} +t$ have Galois group $M_{24}$, as Magma confirms.
We know that the minimal $2$-polynomial of the first polynomial has 2-degree 12, whereas the
minimal $2$-polynomial of the second polynomial has 2-degree 24 (by uniseriality).

\section{Affine Groups}

\noindent A similar result can be obtained for affine groups. 
By way of an example, suppose you started with the
 irreducible polynomial $x^8+x^2+tx+1 \in \F_2(t)[x]$.
Its Galois group is the affine group $AGL(3,2)$ of order 1344, which acts
3-transitively on 8 points.
How would you then construct the binary Hamming code of length 8?

First construct the minimal linearization of $f=x^8+x^2+tx+1$, which is 
\[
L=x^{16}+x^8 + x^4 + (t^2+ 1)x^2 + tx.
\]
There is a module homomorphism from the abstract permutation module $W$ of dimension 8
to the root space $V$ of $L$.
The kernel is a submodule of $W$, of dimension 4.
Since the Galois group is 3-transitive, it follows that  fewer than four roots cannot sum to 0.
It is then straightforward to see that the kernel is the Hamming code.
We omit the details, the argument goes along similar lines to that for $M_{24}$.

A paper of Sin shows that the permutation module for $AGL(n,q)$ is uniserial.
Similar arguments can be made in the general case.

\section{A Projective Polynomial and $PSL(2,q)$ }

\noindent The next lemma provides a family of examples where $m(f)=\deg (f) -1$.
Many of these polynomials have Galois group $PSL(2,p)$.

\begin{lemma}\label{divp1}
Let $p$ be a prime number and let 
 $F$ be a field of characteristic $p$.
Let $f=x^{p+1}-ax-b$ be an irreducible element of $F[x]$.
Then $m(f)=p$.
\end{lemma}

\begin{proof}
Let $\alpha$ be a root of $f(T)$ in some extension field of $F$.
We use the basis $\alpha^k$, for $0\le k \le p$, for the field extension $F(\alpha)$.

We will show that the
elements $\alpha^{p^i}$ for $0\le i \le p$ are linearly dependent over $F$. 
To show this, we will show that these $p+1$ elements lie in a $p$-dimensional 
subspace of the field extension $F(\alpha)$.
To be precise, this $p$-dimensional subspace is the span of 
the basis elements $\alpha^k$, for $1\le k \le p$,
In other words, we claim that when we express each $\alpha^{p^i}$ for $0\le i \le p$ 
as a linear combination of the basis elements $\alpha^k$,  $0\le k \le p$, 
the basis element 1 is not used.
This suffices to prove the theorem.

Since $\alpha^{p+1}=A\alpha+B$ we have
\begin{equation}\label{starter}
\alpha^{i(p+1)}=\sum_{r=0}^i  \binom{i}{r} A^r B^{i-r} \alpha^{r}
\end{equation}
for any $i$ with $1 \le i \le p-1$.
This writes $\alpha^{i(p+1)}$  as a linear combination of $\alpha^k$, $0\le k \le i$.
In particular choosing $i=p-1$ gives
\[
\alpha^{p^2-1}=\sum_{r=0}^{p-1}  \binom{p-1}{r} A^r B^{p-1-r} \alpha^{r}
\]
which implies
\[
\alpha^{p^2}=\sum_{r=0}^{p-1}  \binom{p-1}{r} A^r B^{p-1-r} \alpha^{r+1}.
\]
This writes $\alpha^{p^2}$ as a linear combination of $\alpha^k$, $1\le k \le p$.
Note that the basis element 1 is not used.

To write $\alpha^{p^3}$ as a linear combination of the basis elements,
we raise the $\alpha^{p^2}$ linear combination to the power of $p$.
We will get a term  $\alpha^{ip}$ for each $i$ with $1\le i \le p$.
The terms with  $i=p$ is already taken care of, so it suffices
to show how each term $\alpha^{ip}$  with $1< i < p$ is written in terms of the basis.
We use \eqref{starter} to write
\[
\alpha^{ip}=\alpha^{(i-1)(p+1)+p-i+1}=\alpha^{(i-1)(p+1)}\cdot \alpha^{p-(i-1)}=
\sum_{r=0}^{i-1}  \binom{i-1}{r} A^r B^{i-1-r} \alpha^{r+p-(i-1)}
\]
which expresses $\alpha^{ip}$ as a linear combination of $\alpha^k$, $p-i+1\le k \le p$.
Putting these together, this writes $\alpha^{p^3}$ as a linear combination of $\alpha^k$, $1\le k \le p$.
Note that the basis element 1 is not used.

Inductively, the same argument shows that $\alpha^{p^i}$, $i\ge 4$, is a 
linear combination of $\alpha^k$, $1\le k \le p$.

This shows that $m(f)\le p$. To show equality, note that all the coefficients
of the basis elements are nonzero, so all basis elements $\alpha^k$, $1\le k \le p$ are needed.
\end{proof}

 {\bf Example:} 
 It is shown in Abhyankar \cite{Ab} (in an appendix by Serre) that 
 the Galois group of $f=x^{p+1}+tx+1\in \F_p(t)[x]$ is $PSL(2,p)$.
 In this case $\ker (\psi)$ has dimension 1 and the dimension of $V$ is $p$, by Lemma \ref{divp1}.
 In fact $V$ is the Steinberg module for $PSL(2,p)$.
Thus we have a realization of the Steinberg module (which is irreducible) for $PSL(2,p)$ as
the root space of $f$.


\begin{thebibliography}{99}

\bibitem{Ab}  Abhyankar, S. S.,
Galois theory on the line in nonzero characteristic, Bull. Amer. Math. Soc. {\bf 27} (1992), 68-133.



\bibitem{AY} Abhyankar, S. S.,  Yie, I., 
Some more Mathieu group coverings in characteristic two,
Proc. Amer. Math. Soc {\bf 123} (1994), 1007-1014.



\bibitem{AM} Assmus, E. F., Mattson, H. F., Perfect codes and the Mathieu groups,
Arch. Math. {\bf 17} (1966), 122-135.  

\bibitem{BSEM}
Bary-Soroker, L.,   Entin, A.,  McKemmie, E.,
Galois groups of random additive polynomials,
   Trans.  Amer. Math. Soc. {\bf 377} (2024), 2231-2259. 


\bibitem{BDEPS} Berry, N.,  Dubickas, A.,  Elkies, N.,  Poonen, B.,  Smyth, C.,
The conjugate dimension of algebraic numbers, Quart. J. Math. {\bf 55} (2004), 237-252.


\bibitem{C1} Carmichael, R. D.,
Tactical configurations of rank two,
Amer. J. Math. {\bf 53} (1931), 217-240.

\bibitem{CMT} Conway, J.,  McKay, J.,  Trojan, A.,
Galois groups over function fields of positive characteristic,
Proc. Amer. Math. Soc.
{\bf 138} (2010), 1205-1212.





\bibitem{EntinPopov} Entin, A., Popov, A., 
Probabilistic Galois theory in function fields, 
      arXiv 2311.14862



\bibitem{Gi} Girstmair, K., Linear relations between roots of polynomials, 
Acta Arith. {\bf 89} (1999), 53-96.

\bibitem{Goss} Goss, D., Basic structures of function field arithmetic, 
Springer 1997.




\bibitem{I} Ivanov, A. A., The monster group and Majorana involutions,
Cambridge University Press, 2009.




\bibitem{LN} Lidl, R, Niederreiter, H., Finite fields. Addison-Wesley, Reading, Massachusetts, 1983.

\bibitem{MS} MacWilliams, F. J.,  Sloane, N. J. A., 
The theory of error-correcting codes, North-Holland Mathematical Library, Volume 16, 1977.

\bibitem{P} Paige, L. J., A note on the Mathieu groups, Canad. J. Math. {\bf 9} (1957), 15-18.



\bibitem{T1} Todd, J. A., 
On representations of the Mathieu groups as collineation groups,
J. London Math. Soc. {\bf 34} (1959), 406-416. 

\bibitem{T2} Todd, J. A., 
A representation of the Mathieu group $M_{24}$  as a collineation group,
Ann. di Math. Pure ed Appl. {\bf 71} (1966), 199-238. 



\end{thebibliography}
\end{document}